\documentclass[11pt]{article}

\usepackage{hyperref}
\usepackage[latin1,utf8x]{inputenc}
\usepackage[english]{babel}
\usepackage{amssymb,amsmath,amsfonts,amsthm}

\theoremstyle{plain}
\newtheorem{teo}{Theorem}

\newtheorem{cor}{Corollary}
\newtheorem{lem}{Lemma}
\newtheorem{proposizione}{Proposition}
\newtheorem*{teo*}{Theorem}
\theoremstyle{definition}
\newtheorem{definizione}{Definition}
\newtheorem*{definizionenonum}{Definition}
\newtheorem{oss}{Remark}[section]
\newtheorem{prop}{Proposition}

\newtheorem{ex}{Example}[section]
\newtheorem{ese}{Exercise}[section]
\newtheorem{conj}{Conjecture}

\newcommand{\virgolette}[1]{``#1''}
\newcommand{\bpro}{\begin{proposizione}}		
\newcommand{\epro}{\end{proposizione}}		
\newcommand{\bdef}{\begin{definizione}}		
\newcommand{\bteo}{\begin{teo}}			
\newcommand{\eteo}{\end{teo}}			
\newcommand{\bconj}{\begin{conj}}			
\newcommand{\econj}{\end{conj}}			
\newcommand{\boss}{\begin{oss}}			
\newcommand{\eoss}{\end{oss}}			
\newcommand{\bprop}{\begin{prop}}                        
\newcommand{\eprop}{\end{prop}}                           
\newcommand{\bcor}{\begin{cor}}			
\newcommand{\ecor}{\end{cor}}			
\newcommand{\bdimo}{\begin{proof}}			
\newcommand{\edimo}{\end{proof}}			
\newcommand{\bequ}{\begin{equation}}		
\newcommand{\eequ}{\end{equation}}			
\newcommand{\blem}{\begin{lem}}			
\newcommand{\elem}{\end{lem}}			
\newcommand{\bex}{\begin{ex}}			
\newcommand{\eex}{\end{ex}}				
\newcommand{\bese}{\begin{ese}}			
\newcommand{\eese}{\end{ese}}			
\newcommand{\N}{\mathbb{N}}				
\newcommand{\Z}{\mathbb{Z}}				
\newcommand{\z}{\textbf{\z}}				

\begin{document}
\hsize=14cm


\author{Francesco Monopoli \vspace{0.4cm} \\
\small Dip. di Matematica \virgolette{Federigo Enriques} \\
\small Università degli Studi di Milano \\
\small Via Saldini 50, I-20133 Milano, Italy \\
\small francesco.monopoli@unimi.it
}
\title{Absolute differences \\ along Hamiltonian paths }
\date{}
\maketitle


\begin{abstract}
Given a set $A$ of real numbers consider the complete graph on the elements of $A$. We prove that if $A$ is an arithmetic progression then for every vertex $a\in A$ there exists an hamiltonian path such that the absolute differences of consecutive vertices are pairwise distinct. This result partially proves a conjecture by Zhi-Wei Sun.
\end{abstract}

\renewcommand\refname{References}

\maketitle

\begin{section}*{}
In this paper we consider the following conjecture posed by Z.-W. Sun, formulated among other open problems in \cite[Conjecture 3.1]{zwsun}.

\bconj\label{daconj} Let $A=\{a_1, a_2, \dots, a_n\}$ be a set of $n$ distinct real numbers. Then there is a permutation $b_1, b_2, \dots, b_n$ of $a_1, \dots, a_n$ with $b_1 = a_1$ such that the $n-1$ numbers
$$|b_2 - b_1|, |b_3 - b_2|, \dots, |b_{n} - b_{n-1}|$$
are pairwise distinct.
\econj

Similar problems have been studied in \cite{lev} and \cite{rosa}.

Considering the complete graph on $\{ a_1, a_2, \dots, a_n\}$ and color the edges so that two edges $a_ia_j$ and $a_ka_l$ have the same colour if and only if $|a_i - a_j| = |a_k - a_l|$, in order for the conjecture to be true we need to find for every element $a_h$ a totally multicoloured hamiltonian path starting at $a_h$.

As Z.-W. Sun already pointed out in \cite[Theorem 1.1]{zwsun}, ordering the elements $a_1 < a_2 < \dots < a_n$ we can easily find such an hamiltonian path starting from $a_1$ or $a_n$: if $n =2k$ is even we can consider the permutation $$(a_1, a_n, a_2, a_{n-1}, \dots, a_{k-1}, a_{k+2}, a_k, a_{k+1}),$$ and if $n=2k-1$ is odd consider the permutation $$(a_1, a_n, a_2, a_{n-1}, \dots, a_{k-1}, a_{k+1}, a_k).$$

If the cardinality of $A-A = \{ a_i - a_j : a_i, a_j \in A\}$ is large, then heuristically it should be easy to find hamiltonian paths as require, whereas this should be an harder task for structured sets, where $|A-A|$ can be as small as $|A|-1$.

However, we are able to prove that the conjecture holds in these cases.

\bteo\label{dateo}
Conjecture \ref{daconj} holds if $A$ is an arithmetic progression.
\eteo

Moreover, as expected, the conjecture holds if the set $A$ does not have a particular additive structure.

\bteo\label{dateo2}
Conjecture \ref{daconj} holds if $E(A, A) = c|A|^2$ for $c < 5/2$, where $E(A,A)$ is the additive energy of the set $A$.
\eteo
\end{section}

\begin{section}*{Arithmetic progressions}
In order to prove \ref{dateo} let without loss of generality $A=[n]:=[1, n] \cap \Z$ be the set of the first $n$ positive integers.

Fix an element $s \in [n]$. We want to find a permutation $a=(a_1, \dots, a_n)$ of $[1,n]$ with $a_1 = s$ such that the $n-1$ differences
$$|a_2-a_1|, \dots, |a_{n}-a_{n-1}|$$
are pairwise distinct.

Define the set of absolute differences of the sequence $a$ as $d(a) := \{ |a_{i+1}-a_{i}|: a_i, a_{i+1} \in a, a_i \neq a_j\}$. We want to find an $a$ of $[1,n]$ such that $|d(a)|=n-1$.

\begin{definizionenonum}
We call a permutation $a=(a_1, \dots, a_n)$ of $[n]$ a \emph{good sequence} if the differences $|a_2-a_1|, \dots, |a_n-a_{n-1}|$ are pairwise distinct and one of the following holds:
\begin{enumerate}

\item $a_{2l+1} \geq \left\lceil \frac{n+1}{2}\right\rceil$ and $a_{2l} < \left\lceil\frac{n+1}{2} \right\rceil$ whenever $ 2l+1, 2l \in  [n]$.
\item $a_{2l+1} \leq \left\lfloor \frac{n+1}{2} \right\rfloor$ and $a_{2l} > \left\lfloor \frac{n+1}{2} \right\rfloor$ whenever $2l+1, 2l \in [n]$.
\end{enumerate}

\end{definizionenonum}

Clearly, if we can find such a sequence, then theorem \ref{dateo} would be proved. Unfortunately, for some starting points this is not possible, but we will be able to treat those separately.

These kind of special  permutations are useful because they allow us to build new good sequences with different starting points, with the two procedures explained in the following lemma.

\blem\label{lem1}
Let $a=(a_1, \dots, a_n)$ be a good sequence. Then the following hold:
\begin{enumerate}
\item The sequence $b=\{b_i \}$ given by $b_i = n+1-a_i$ is again a good sequence

\item Suppose $a_{2l+1} \geq \left\lceil \frac{n+1}{2} \right\rceil$. 
Then the sequence $b=\{b_i \}$ given by $b_{2l+1} = a_{2l+1} - \left\lfloor \frac{n}{2} \right\rfloor, b_{2l} =a_{2l} + \left\lfloor \frac{n+1}{2} \right\rfloor$ is again a good sequence.
\item Suppose $a_{2l+1} \leq \left\lfloor \frac{n+1}{2} \right\rfloor$. 
Then the sequence $b=\{b_i \}$ given by $b_{2l+1} = a_{2l+1} + \left\lfloor \frac{n}{2} \right\rfloor, b_{2l} =a_{2l} - \left\lfloor \frac{n+1}{2} \right\rfloor$ is again a good sequence.
\end{enumerate}
\elem

\bdimo

\begin{enumerate}
\item Suppose $a_{2l+1} \leq \left\lfloor \frac{n+1}{2} \right\rfloor$. Then $b_{2l+1} = n+1-a_{2l+1} \geq n+1 - \left\lfloor \frac{n+1}{2} \right\rfloor = \left\lceil \frac{n+1}{2} \right\rceil$, and $b_{2l} = n+1-a_{2l} <\left\lceil \frac{n+1}{2} \right\rceil$.

A similar statement holds if $a_{2l+1} \geq \left\lfloor \frac{n+1}{2} \right\rfloor$.

\item
Since $\left\lceil \frac{n+1}{2} \right\rceil \leq a_{2l+1} \leq n$ and $0 < a_{2l} < \left\lceil \frac{n+1}{2} \right\rceil$, we have 
$$1 = \left\lceil\frac{n+1}{2} \right\rceil - \left\lfloor\frac{n}{2} \right\rfloor \leq b_{2l+1} \leq n- \left\lfloor \frac{n}{2} \right\rfloor = \left\lfloor \frac{n+1}{2} \right\rfloor,$$

$$ \left\lfloor \frac{n+1}{2} \right\rfloor < b_{2l} < \left\lfloor \frac{n+1}{2} \right\rfloor + \left\lceil \frac{n+1}{2} \right\rceil = n+1.$$

This shows that the new sequence $b$ is again a permutation of $[1,n]$ and satisfies the second condition for being a good sequence.

Moreover, the differences $|b_2-b_1|, \dots, |b_n-b_{n-1}|$ are pairwise disjoint: 
\begin{eqnarray*}
|b_{2l+1} - b_{2l}|  &=& |a_{2l+1} - a_{2l} - \left\lfloor \frac{n}{2} \right\rfloor - \left\lfloor \frac{n+1}{2} \right\rfloor| \\
			&=&n-|a_{2l+1}-a_{2l}|,
\end{eqnarray*}
and so the differences $|b_2-b_1|, \dots, |b_n-b_{n-1}|$ are just a permutation of the elements $\{|a_{i+1}-a_{i}|\}$, which were pairwise disjoint by hypothesis, and so $b$ is again a good sequence.
\item Same proof as in point 2.
\end{enumerate}
\edimo

We can now prove the main result, which clearly implies theorem \ref{dateo}.
\bteo\label{teo1}
If $n \not\equiv 1$ mod $4$ then for every $s \in [n]$ there exists a good permutation $a=(a_1, \dots, a_n)$ of $[n]$ with $a_1=s$.

If $n \equiv 1$ mod $4$ then for every $s \in [n]$ there exists a permutation $a=(a_1, \dots, a_n)$ of $[n]$ with $a_1 =s$ and $|d(a)|=n-1$. Moreover, if $s \neq \frac{1}{2}\left( \left\lfloor \frac{n+1}{2}\right\rfloor+1 \right)$ one can find a good sequence starting from $s$.
\eteo

\bdimo
The proof goes by induction on $n$. Because of the first part of lemma \ref{lem1} we can prove it just for starting points $s \leq \left\lfloor \frac{n+1}{2}\right\rfloor$.

If $s=1$ the sequence $(1, n, 2, n-1, \dots, \left\lfloor \frac{n+1}{2}\right\rfloor + \delta)$, where $\delta = 1$ if $n$ is even and $\delta = 0 $ if $n$ is odd, is clearly a good sequence.

Take $2 \leq s \leq \frac{1}{2} \left\lfloor \frac{n+1}{2}\right\rfloor$. We consider two cases:

Case 1: $n-2s \not\equiv 1$ mod $4$.

Then we can consider the following sequence:
$$b=(s, n-s+1, s-1, n-s+2, \dots, 1, n),$$
with $d(b) = [n-2s+1, n-1]$.

We choose the next element as $\alpha = 2s \leq \left\lfloor \frac{n+1}{2}\right\rfloor$ in order to get the absolute difference $n-2s$.

By induction hypothesis we can find a good permutation $c$ of $[1,n-2s]$ starting from $s$, so that $d(c)=[1,n-2s-1]$. Since
$$ s \leq \frac{1}{2}\left\lfloor \frac{n+1}{2}\right\rfloor = \frac{1}{2}\left\lfloor \frac{n-2s+1}{2}\right\rfloor + \frac{s}{2}$$
implies $$s \leq \left\lfloor \frac{n-2s+1}{2}\right\rfloor$$ we have that the permutation $a$ obtained by linking together $b$ and $s+c=(c_1 + s=2s, c_2+s, \dots )$ satisfies
$$a_{2l+1} \leq \left\lfloor \frac{n+1}{2}\right\rfloor, \qquad a_{2l} > \left\lfloor \frac{n+1}{2}\right\rfloor,$$
so that $a$ is a good sequence starting from $s$.

Case 2: $n-2s+1 \not\equiv 1$ mod $4$.

In this case we start from the sequence

$$b=(s, n-s+2, s-1, n-s+3, \dots, 2, n, 1),$$

so that $d(b) = [n-2s+2, n-1]$, and take the next element as $\alpha=n-2s+2$ in order to get the absolute difference $n-2s+1$.

Using the inductive hypothesis we find a good permutation $c$ of $[1, n-2s+1]$ starting from $n-3s+2$, with $d(c)=[1,n-2s]$.

Since our hypothesis on $s$ imply that $n-3s+2 \geq \left\lceil \frac{n-2s+2}{2}\right\rceil $ we get that 
$$s+c_{2l+1} \geq \left\lceil \frac{n+2}{2}\right\rceil = \left\lfloor \frac{n+1}{2}\right\rfloor  +1 \qquad
s+c_{2l} < \left\lfloor \frac{n+1}{2}\right\rfloor  +1,$$

so that the sequence $a$ obtained by chaining $b$ and $s+c$ is indeed a good sequence.

Since for every $n$ and $s \leq \frac{1}{2}\left\lfloor \frac{n+1}{2}\right\rfloor$ either $n-2s$ or $n-2s+1$ is not congruent to $1$ modulo $4$, the result is proven in these cases.

Suppose now $\frac{1}{2}\left\lfloor \frac{n+1}{2}\right\rfloor  < s \leq  \left\lfloor \frac{n+1}{2}\right\rfloor$. Then $s' = n+1-s- \left\lfloor \frac{n}{2}\right\rfloor <  \frac{1}{2} \left\lfloor \frac{n+1}{2}\right\rfloor+1$, and by lemma \ref{lem1} we are done unless $n \equiv 1,2$ modulo $4$ and $s=\frac{1}{2} \left\lfloor \frac{n+1}{2}\right\rfloor + \frac{1}{2}.$

We study these cases separately.

Case 1: $n \equiv 1$ mod $4$.

We consider the sequence
$$(s, n-s+2, s-1, n-s+3, \dots, 2, n, 1, n-2s+2= \left\lceil \frac{n+1}{2}\right\rceil , $$ $$\left\lceil \frac{n+1}{2}\right\rceil +1, \left\lceil \frac{n+1}{2}\right\rceil -1, \dots,n- s+1).$$

Case 2: $n \equiv 2$ mod $4$.

We consider the sequence
$$(s, n-s+2, s-1, n-s+3, \dots, 2, n, 1, n-2s+2= \left\lceil \frac{n+1}{2}\right\rceil ,$$ $$ \left\lceil \frac{n+1}{2}\right\rceil -1, \left\lceil \frac{n+1}{2}\right\rceil +1, \dots, s+1, n-s+1).$$

\edimo

Given a hamiltonian path starting from $s$, something more can be said about the parity of the ending point of such a permutation.

\bcor
$n$ is congruent to $0$ or $1$ modulo $4$ if and only if for every permutation $(a_1, \dots, a_n)$ of $[n]$ such that the $n-1$ differences $|a_2-a_1|, \dots, |a_n - a_{n-1}|$ are different, we have $a_1 \equiv a_n$ modulo $2$.
\ecor
\bdimo
Given a permutation of $[n]$ with the property described in the statement, whose existence is guaranteed by theorem \ref{teo1} for every starting point $a_1$, we have
$$\sum_{i=2}^n |a_i - a_{i-1}| \equiv \sum_{i=2}^n a_i - a_{i-1} \equiv a_n - a_1 \qquad \mbox{mod $2$}.$$
On the other hand, the LHS is equal to $\sum_{i=1}^n i = \frac{n(n-1)}{2},$
which is congruent to $0$ modulo $2$ if and only if $n$ is congruent to $0$ or $1$ modulo $4$.
\edimo

\end{section}

\begin{section}*{Random sets}
In this section we show that the conjecture holds for ``random" sets, i.e. sets whose additive energy $E(A, A)$ is small.

Consider the similar problem of finding a hamiltonian cycle  $a=(a_1, \dots, a_n)$ of the elements of a set $A$ of cardinality $n$ such that the $n$ differences $|a_2 - a_1|, \dots, |a_n - a_{n-1}|, |a_1 - a_n|$ are pairwise disjoint. If we were able to find such a permutation then clearly we would have found, for any starting point $s \in A$ a hamiltonian path satisfying our original condition.

Of course this is not always the case, since for example if $A$ is an arithmetic progression o length $n$, such a cycle cannot exist, for $|A-A|_+ := |(A-A) \cap \N| = n-1$ and therefore it's impossible to produce $n$ distinct absolute differences.

Moreover, even sets $A$, $|A|=n$, with $|A-A|_+ \geq |A|$ might fail to satisfy this condition: let $A \subseteq [n+1], |A-A|_+ = n$ and suppose $a=(a_1, \dots, a_n)$ is an hamiltonian path on $A$ with distinct consecutive absolute differences. Then
$$
\frac{n(n+1)}{2} = \sum_{i=1}^n i =\sum_{i=2}^{n} |a_i - a_{i-1}| + |a_1 - a_n| \equiv \sum_{i=2}^{n} a_i - a_{i-1} + a_1 - a_n \equiv 0
$$
modulo $2$, and if $n$ is congruent to $1$ or $2$ modulo $4$ this cannot happen. 

Given a random circular permutation $a=(a_1, \dots, a_n)$ of $A$ let $d(a)=\{ |a_2 - a_1|, \dots, |a_n - a_{n-1}|, |a_1 - a_n| \}.$

Then
\begin{equation}\label{eq1} E(|d(a)|) = \sum_{d \in (A-A)_+} P(d \in d(a)).\end{equation}

Fix $d \in (A-A)_+$. Let $X_i$ be the event $|a_i - a_{i-1}|=d$ for $i=2, \dots, n$, and $X_1$ be the event $|a_1 - a_n| = d$. 

Then
\begin{equation}\label{eq2}P(d \in d(a)) = P(X_1 \vee \dots \vee X_n) \geq \sum_{i=1}^n P(X_i) - \sum_{1\leq i < j \leq n} P(X_i\wedge X_j)\end{equation}
by inclusion-exclusion.

Let $s(d) = |\{ a \in A : a-d, a+d \in A \}|$. be the number of $3$-terms arithmetic progressions of difference $d$ contained in $A$, and $r_{A, -A}(x):= | \{ (a, a') \subseteq A \times A : a-a' = x\} |$ An elementary estimate is the following

\blem
Let $|A|=n$. Then $\sum_{d\in |A-A|_+} s(d) \leq n^2/4$.
\elem
\bdimo
If $A=\{a_1 < \dots <a_n\}$ then $a_i$ can be the middle term of no more than $\min(i-1, n-i)$ three terms arithmetic progressions. Hence
$$\sum_{d\in |A-A|_+} s(d)  \leq 2\sum_{i=1}^{\left\lfloor \frac{n+1}{2}\right\rfloor} (i-1) \leq \frac{n^2}{4} $$
\edimo

Then
\begin{eqnarray*}
P(X_i) &=&\frac{2r_{A, -A}(d)}{n(n-1)}\\
P(X_i \wedge X_{i+1}) &=& \frac{2s(d)}{n(n-1)(n-2)} \\
P(X_i \wedge X_j, i+1<j) &=& \frac{4r_{A, -A}(d)(r_{A, -A}(d)-1)}{n(n-1)(n-2)(n-3)} 
\end{eqnarray*}

Putting these equalities in \ref{eq1} and \ref{eq2} we get

\begin{eqnarray*}
E(|d(a)|) &\geq&\frac{2\sum_{d\in (A-A)_+} r_{A, -A}(d)}{n-1} - \frac{2\sum_{d\in (A-A)_+}s(d)}{(n-1)(n-2)} + \\
	& &- \frac{2\sum_{d \in (A-A)_+} r_{A, -A}(d) (r_{A, -A}(d) -1)}{n(n-3)} \\
	&\geq& n - \frac{1}{2}\frac{n^2}{(n-1)(n-2)} - \frac{E(A, A) - 2n^2 +n}{n(n-3)}.
\end{eqnarray*}

Then, for $n\gg 1$ and a set $A$ with $E(A, A) =c n^2$ for a $c < 5/2$, we have $E(|d(a)|) > n-1$, and hence there exists a hamiltonian cycle $a$ of $A$ with $|d(a)|=n$ as required.
\end{section}

\end{document}